\documentclass[12pt]{article}
\usepackage{amssymb}
\usepackage{amsmath}

\begin{document}

\begin{center}
\textbf{\large Tame semicascades and cascades generated\\ by
affine self-mappings of the $d$-torus}
\end{center}

\begin{center}
\textsc{Vladimir Lebedev}
\end{center}

\begin{quotation}
{\small \textbf{Abstract.} We give a complete characterization
of the affine self-mappings $\varphi$ of the torus $\mathbb
T^d$ that generate K\"ohler-tame semicascades and cascades.
Namely, we show that the semicascade generated by $\varphi$ is
tame if and only if the matrix $A$ of $\varphi$ satisfies
$A^p=A^q$, where $p$ and $q$ are some nonnegative integers,
$p\neq q$. For cascades the corresponding condition has the
form $A^m=I$, where $m$ is some positive integer and $I$ is the
identity matrix.


}
\end{quotation}

\quad

\textbf{1. Introduction and statement of results.} A sequence
$f_n, \,n=1, 2, \ldots,$ of real-valued functions on a set $X$
is said to be independent if there exist $a, b\in\mathbb R,
\,a<b,$ such that
$$
\bigcap_{n\in P} \{x : f_n(x)<a\}\cap
\bigcap_{n\in Q} \{x : f_n(x)>b\}\neq\varnothing
$$
for all finite disjoint subsets $P, Q$ of indices.

Let $X$ be a compact Hausdorff space and let $G$ be a certain
group or a semigroup whose elements are continuous mappings of
$X$ into itself. The dynamical system $(X, G)$ is called
\emph{tame} if for every real-valued continuous function $f$ on
$X$ the family $\{f_g: \,g\in G\}$ does not contain an
independent sequence [2,~Def.~3.1]. Here $f_g$ stands for the
function $f_g(x)=f(gx), \,x\in X$.

This definition of a tame (originally called regular) dynamical
system goes back to K\"ohler, who considered semicascades and
obtained first results on tameness [3, Sec. 5]. The study
continued by Glasner, Megrelishvili, and other authors, see the
recent survey [2] and the references therein, showed that the
tame/untame dichotomy is closely connected with such properties
of a system as minimality, distality, nonsensitivity, almost
periodicity, chaotic behaviour, certain ergodic properties, and
the structure of the Ellis enveloping semigroup. We note that
it is usually difficult to decide whether a given system is
tame.

Several conditions are known to be equivalent to tameness. In
the important case when $X$ is a compact metric space the above
K\"ohler definition of a tame system reduces to a very simple
equivalent definition below which reflects the structure of
orbits $\mathcal{O}(x)=\{gx: \,g\in G\}, \,x\in X,$ of a
system, and is especially meaningful for semicascades and
cascades. (We provide a short explanation of the equivalence in
the concluding Section 3, see Remark 1.)

\quad

\textsc{Definition} (an equivalent version in the metric case).
Let $X$ be a compact metric space. We say that a system $(X,
G)$ is tame if each sequence $g_n\in G, \,n=1, 2, \ldots,$ has
a pointwise convergent subsequence $g_{n_k}, \,k=1, 2, \ldots$
(i.e., a subsequence such that for every $x\in X$ the sequence
$g_{n_1}x, g_{n_2}x, \ldots$ converges in $X$).

\quad

Let $\varphi : X\rightarrow X$ be a continuous mapping. We
recall that the semicascade generated by $\varphi$ is the
system $(X, G_\varphi^+)$, where $G_\varphi^+=\{\varphi^n,
\,n=0, 1, 2, \ldots\}$ is the semigroup of iterations of
$\varphi$. Here $\varphi^0$ is the identity mapping and
$\varphi^{n+1}=\varphi\circ\varphi^n, \,n=0, 1, 2, \ldots$.
When $\varphi$ is a self-homeomorphism of $X$ we can consider
the cascade generated by $\varphi$, i.e., the system $(X,
G_\varphi)$, where $G_\varphi=\{\varphi^n, n\in\mathbb Z\}$ is
the group of iterations of $\varphi$ (naturally $\varphi^n$ for
$n<0$ are the iterations of the inverse mapping
$\varphi^{-1}$).

In this note we consider semicascades and cascades generated by
affine self-mappings $\varphi$ of the torus $\mathbb
T^d=\mathbb R^d/(2\pi\mathbb Z)^d$ (as usual $\mathbb R$
denotes the real line and $\mathbb Z$ is its additive subgroup
of integers.). Given a $\varphi$, we have $\varphi(x)=Ax+b$,
where $A$ is an integer $d\times d$ matrix and $b\in\mathbb
T^d$. Certainly one can consider a cascade generated by
$\varphi$ only in the case when $A$ is nonsingular and its
inverse $A^{-1}$ is integer, i.e., when $\det A=\pm1$.

How to distinguish if the semicascade (cascade) generated by an
affine self-mapping of the torus is tame? This question was
suggested by A. V. Romanov (private communication). The results
of this work are the following two theorems.

\quad

\textsc{Theorem 1.} \emph{The semicascade generated by an
affine self-mapping $\varphi$ of $\mathbb T^d$ is tame if and
only if the matrix $A$ of $\varphi$ satisfies the condition
$A^p=A^q$, where $p$ and $q$ are some nonnegative integers with
$p\neq q$.}

\quad

\textsc{Theorem 2.} \emph{The cascade generated by an affine
self-mapping $\varphi$ of $\mathbb T^d$ is tame if and only if
the matrix $A$ of $\varphi$ satisfies the condition $A^m=I$,
where $m$ is some positive integer and $I$ is the identity
matrix.}

\quad

The proof is given in the next section. It is based on
elementary Fourier analysis argument.

The results of this work were announced at the conference
``Topology, Geometry, and Dynamics: Rokhlin – 100'' [4].

\quad

\textbf{2. Proof of the theorems.} It is convenient to consider
the general case of families (not necessarily semigroups or
groups) of affine self-mappings of the torus. Theorems 1 and 2
immediately follow from the lemma below.

\quad

\textsc{Lemma.} \emph{Let $\Phi$ be a family of affine
self-mappings of $\mathbb T^d$. Let $M(\Phi)$ be the family of
all matrices of the mappings in $\Phi$. The following
conditions are equivalent:\\
\emph{(i)} each sequence in $\Phi$ has a pointwise
convergent subsequence;\\
\emph{(ii)} the family $M(\Phi)$ is finite.}

\quad

\textsc{Proof of the Lemma.} The part $\mathrm{(ii)\Rightarrow
(i)}$ is trivial. Let us show that $\mathrm{(i)\Rightarrow
(ii)}$. Given a Lebesgue integrable function $f$ on $\mathbb
T^d$ its Fourier transform $\widehat{f}$ is defined by
$$
\widehat{f}(\lambda)=\frac{1}{(2\pi)^d}\int_{\mathbb T^d}f(x)e^{-i(\lambda, \,x)} dx,
\qquad \lambda\in\mathbb Z^d,
$$
where $(\lambda, x)$ is the usual inner product of vectors
$\lambda\in\mathbb Z^d$ and $x\in\mathbb T^d$, i.e., $(\lambda,
x)=\sum_{j=1}^d\lambda^j x^j$ for $\lambda=(\lambda^1,
\lambda^2, \ldots, \lambda^d), \, x=(x^1, x^2, \ldots, x^d)$.

For a vector $\lambda\in\mathbb Z^d$ let $e_\lambda$ stands for
the exponential function with frequency $\lambda$, i.e.,
$e_\lambda(x)=e^{i(\lambda, \,x)}, \, x\in\mathbb T^d$. Note,
that if $\gamma_1, \gamma_2, \ldots$ is a sequence of real
numbers, and $\lambda_1, \lambda_2, \ldots$ is an unbounded
sequence of vectors in $\mathbb Z^d$, then the sequence
$e^{i\gamma_1}e_{\lambda_1}, e^{i\gamma_2}e_{\lambda_2},
\ldots$ is not pointwise convergent on $\mathbb T^d$. Indeed,
assuming the contrary, let $\lim_{n\rightarrow\infty}
e^{i\gamma_n}e^{i(\lambda_n, x)}=\xi(x), \,x\in\mathbb T^d$.
Using dominated convergence theorem we see that the function
$\xi$ is integrable on $\mathbb T^d$ and, since
$\xi(x)e^{-i\gamma_n}e^{-i(\lambda_n, x)}\rightarrow 1$ for all
$x\in\mathbb T^d$, we obtain
$$
\frac{1}{(2\pi)^d}\int_{\mathbb
T^d}\xi(x)e^{-i\gamma_n}e^{-i(\lambda_n, \,x)} dx\rightarrow 1,
\quad\mathrm{as} \,\, n\rightarrow\infty,
$$
so, $e^{-i\gamma_n}\widehat{\xi}(\lambda_n)\rightarrow 1$,
which is impossible since the Fourier coefficients of every
integrable function vanish at infinity.

Given a matrix $A$, denote its transpose by $A^*$. Consider an
arbitrary vector $u\in\mathbb Z^d$ and the set
$E(u)=\{e_u\circ\varphi: \,\varphi\in \Phi\}$. Certainly, (i)
implies that each sequence in $E(u)$ has a pointwise convergent
subsequence. At the same time for $\varphi(x)=Ax+b$ we have
$e_u\circ\varphi(x)=e^{i(u, Ax+b)}=e^{i(u, b)}e^{i(u,
Ax)}=e^{i(u, b)}e^{i(A^*u, \,x)}=e^{i(u, b)}e_{A^*u}(x)$. Thus
we see that from (i) it follows that the set
$\Lambda(u)=\{A^*u: \,A\in M(\Phi)\}$, which is the set of
frequencies of functions in $E(u)$, is finite.

Consider the vectors $u_j=(\underbrace{0, 0, \ldots, 0, 1}_j,
0, \ldots, 0)\in\mathbb Z^d$, where $j=1, 2, \ldots, d$.
According to what we have shown, for every $j$ the set
$\Lambda(u_j)=\{A^* u_j: \,A\in M(\Phi)\}$ is finite. Since
$A^*u_j$ is the $j$th column of the matrix $A^*$, i.e., the
$j$th row of the matrix $A$, we see that for every $j, \,j=1,
2, \ldots, d,$ the set of the $j$th rows of the matrices $A\in
M(\Phi)$ is finite. Hence, $M(\Phi)$ is finite. The Lemma is
proved. Theorems 1 and 2 follow.

\quad

\textbf{3. Remarks.} 1. Given a Hausdorff compact space $X$ let
$C(X)$ be the Banach space of all complex-valued continuous
functions on $X$ (with the usual sup-norm). Let $F$ be a
bounded set of real-valued functions in $C(X)$. It is known
that the following two conditions are equivalent:

(a) $F$ does not contain an independent sequence;

(b)  each sequence in $F$ has a pointwise convergent
subsequence.

\noindent The equivalence of these conditions as well as their
equivalence to certain other ones originates from the famous
Rosenthal's $l^1$ theorem [5], [6]. For a thorough discussion
see [1] (especially [1, Theorem 3.11]). This equivalence plays
a significant role in investigations related to tameness, see
[2, Sec. 2 and 3]. Applied to the families $\{f_g: \,g \in
G\}$, where $f\in C(X)$ are real functions, it yields the
corresponding equivalent properties of $(X, G)$.

Consider a family $\Phi$ of self-mappings of $X$. Certainly if
each sequence in $\Phi$ has a pointwise convergent subsequence,
then for every real-valued $f\in C(X)$ each sequence in
$\{f\circ\varphi: \,\varphi\in\Phi\}$ has a pointwise
convergent subsequence. One can easily verify that in the
metric case the converse is also true. Indeed, if $X$ is a
compact metric space, then $C(X)$ is separable, so, we can
chose a countable family of real functions $f_1, f_2, \ldots$
in $C(X)$ that distinguishes between the points of $X$, i.e.,
such that for every $u, v\in X, \,u\neq v,$ there exists $j$
with $f_j(u)\neq f_j(v)$. Let $\{\varphi_n\}$ be a sequence in
$\Phi$. We apply the diagonal process, namely, starting from
the sequence $\{\varphi_n^0\}\overset{def}{=}\{\varphi_n\}$ we
construct inductively a family of sequences $\{\varphi_n^j,
\,n=1, 2, \ldots\}, \,j=1, 2, \ldots,$ so that for each $j$ the
sequence $\{\varphi_n^{j}, \,n=1, 2, \ldots\}$ is a subsequence
of $\{\varphi_n^{j-1}, \,n=1, 2, \ldots\}$ and the sequence
$\{f_j\circ\varphi_n^j, \,n=1, 2, \ldots\}$ is pointwise
convergent. Clearly $\{\varphi_n^{(n)}, \,n=1, 2, \ldots\}$ is
a pointwise convergent subsequence of $\{\varphi_n\}$.

This argument, combined with the equivalence of (a) and (b),
shows that in the metric case the definition of tameness given
in terms of pointwise convergent sequences in $G$ is equivalent
to the original one given in terms of independent sequences
(see Introduction).

2. In relation with Theorems 1 and 2 let us indicate two
conditions in terms of the Jordan form $A_J$ of a matrix $A$
equivalent to the $A^p=A^q$ condition. Clearly, if $A^p=A^q$,
where $p, q$ are certain nonnegative integers, $p\neq q$, then
there exists an integer $s>0$ such that each eigenvalue of $A$
is ether $0$ or is an $s$th root of unity and all Jordan blocks
which correspond to nonzero eigenvalues are degenerate, i.e.,
of size $1\times1$. This in turn implies that there exists a
positive integer $m$ such that $(A_J)^m$ is a diagonal matrix
with only ones and zeroes on the diagonal. Obviously the
converse is also true: every matrix $A$ whose Jordan form
satisfies one of the two conditions above satisfies the
$A^p=A^q$ condition.

\quad

\textbf{Acknowledgements.} The author is grateful to A. V.
Romanov for helpful discussions.

The article was prepared within the framework of the Academic
Fund Program at the National Research University Higher School
of Economics (HSE University) in 2019--2020 (grant
No.19-01-008) and within the framework of the Russian Academic
Excellence Project ``5--100''.

\begin{center}
\textbf{References}
\end{center}

\flushleft
\begin{enumerate}

\item D. van Dulst, \emph{Characterizations of Banach spaces
    not containing $l^1$}, Centrum voor Wiskunde en
    Informatica, Amsterdam, 1989.

\item Eli Glasner, Michael Megrelishvili,  ``More on Tame
    Dynamical Systems'', In: Ferenczi S., Kulaga-Przymus J.,
    Lemanczyk M. (eds) Ergodic Theory and Dynamical Systems
    in their Interactions with Arithmetics and Combinatorics.
    Lecture Notes in Mathematics, vol 2213, pp 351--392,
    Springer, Cham, 2018.

\item  Angela K\"ohler, ``Enveloping semigroups for flows'',
    \emph{Proceedings of the Royal Irish Academy}, 95A:2,
    (1995), 179--191.

\item  Vladimir Lebedev, ``Tame semicascades and cascades
    generated by affine self-mappings of the $d$-torus'',
    Topology, Geometry, and Dynamics: Rokhlin – 100; Aug
    19--23, 2019, Euler Int. Math. Inst., St. Petersburg,
    Abstracts, p. 55.

\item Haskell P. Rosenthal, ``A characterization of Banach
    spaces containing $l^1$\,'', Proc. Nat. Acad. Sci.
    U.S.A., 71:6 (1974), 2411--2413.

\item Haskell P. Rosenthal, ``Some recent discoveries in the
    isomorphic theory of Banach spaces'', Bull. Amer. Math.
    Soc., 84:5 (1978), 803--831.

\end{enumerate}

\quad

\qquad School of Applied Mathematics\\
\qquad National Research University Higher School of Economics\\
\qquad 34 Tallinskaya St.\\
\qquad Moscow, 123458 Russia\\
\qquad E-mail address: \emph{lebedevhome@gmail.com}

\end{document}